\documentclass{amsart}

\usepackage{amsmath,amscd}

\let\ssize\scriptstyle
\newsymbol\diagup 231E
\newsymbol\diagdown 231F
\newif\ifFIRST\newdimen\MAXright\MAXright0pt
\def\sdynkin{\bgroup\eightpoint\dynkin}
\def\endsdynkin{\enddynkin\egroup}
\def\dynkin{\bgroup\FIRSTtrue\hskip.5em\setbox1\hbox{$\diagup$}%
\setbox2\hbox{$\diagdown$}%
\setbox0\hbox to2\wd1{\hrulefill}%
\setbox3\hbox{$\bullet$}%
\setbox4\hbox{$\times$}%
\setbox7\hbox{$\circ$}
\def\whiteroot##1{\ifFIRST\setbox5\hbox{$##1$}\ifdim\wd5>1.3em
\hskip-.5em\hskip.5\wd5\fi\fi\FIRSTfalse
\hskip-.25em\raise1.5\wd3\hbox to0pt{\hss\hskip.45em$
\ssize##1$\hss}\copy7\hskip-.25em\setbox6\hbox{$##1$}
\MAXright\wd6}
\def\root##1{\ifFIRST\setbox5\hbox{$##1$}\ifdim\wd5>1.3em%
\hskip-.5em\hskip.5\wd5\fi\fi\FIRSTfalse%
\hskip-.25em\raise1.5\wd3\hbox to0pt{\hss\hskip.45em$%
\ssize##1$\hss}\copy3\hskip-.25em\setbox6\hbox{$##1$}%
\MAXright\wd6}%
\def\whitedroot##1{\ifFIRST\setbox5\hbox{$##1$}\ifdim\wd5>1.3em
\hskip-.5em\hskip.5\wd5\fi\fi\FIRSTfalse
\hskip-.25em\lower1.8\wd3\hbox to0pt{\hss\hskip.45em$
\ssize##1$\hss}\copy7\hskip-.25em\setbox6\hbox{$##1$}
\MAXright\wd6}%
\def\whiterroot##1{\hskip-.25em\copy7\hbox to0pt{\hskip.3em$\ssize##1$\hss}%
\hskip-.25em\setbox6\hbox{\hskip.6em$##1##1$}%
\MAXright\wd6}%
\def\droot##1{\ifFIRST\setbox5\hbox{$##1$}\ifdim\wd5>1.3em%
\hskip-.5em\hskip.5\wd5\fi\fi\FIRSTfalse%
\hskip-.25em\lower1.8\wd3\hbox to0pt{\hss\hskip.45em$%
\ssize##1$\hss}\copy3\hskip-.25em\setbox6\hbox{$##1$}%
\MAXright\wd6}%
\def\rroot##1{\hskip-.25em\copy3\hbox to0pt{\hskip.3em$\ssize##1$\hss}%
\hskip-.25em\setbox6\hbox{\hskip.6em$##1##1$}%
\MAXright\wd6}%
\def\norroot##1{\hskip-.36em\copy4\hbox to0pt{\hskip.3em$\ssize##1$\hss}%
\hskip-.48em\setbox6\hbox{\hskip.6em$##1##1$}%
\MAXright\wd6}%
\def\noroot##1{\ifFIRST\setbox5\hbox{$##1$}\ifdim\wd5>1.3em%
\hskip-.5em\hskip.5\wd5\fi\fi\FIRSTfalse%
\hskip-.36em\raise1.5\wd3\hbox to0pt{\hss\hskip.6em$%
\ssize##1$\hss}\copy4\hskip-.38em\setbox6\hbox{$##1$}%
\MAXright\wd6}%
\def\nodroot##1{\ifFIRST\setbox5\hbox{$##1$}\ifdim\wd5>1.3em%
\hskip-.5em\hskip.5\wd5\fi\fi\FIRSTfalse%
\hskip-.36em\lower1.8\wd3\hbox to0pt{\hss\hskip.6em$%
\ssize##1$\hss}\copy4\hskip-.38em\setbox6\hbox{$##1$}%
\MAXright\wd6}%
\def\nolink{\hskip\wd0}
\def\link{\raise.22em\copy0}%
\def\llink##1{\raise.32em\copy0\hskip-\wd0%
\raise.12em\copy0\hskip-.5\wd0\hbox to0pt{\hss$##1$\hss}\hskip.5\wd0}%
\def\lllink##1{\raise.22em\copy0\hskip-\wd0\raise.32em\copy0\hskip-\wd0%
\raise.12em\copy0\hskip-.5\wd0\hbox to0pt{\hss$##1$\hss}\hskip.5\wd0}%
\def\rootupright##1{\hbox to0pt{\raise.45em\copy1\hskip-.25em\raise1.3\ht1%
\hbox{\copy3\hskip.3em$\ssize##1$}\hss}%
\setbox6\hbox{\hskip.6em\copy1\copy1$##1##1$}%
\ifdim\MAXright<\wd6\MAXright\wd6\fi}%
\def\whiterootupright##1{\hbox to0pt{\raise.45em\copy1\hskip-.25em\raise1.3\ht1
\hbox{\copy7\hskip.3em$\ssize##1$}\hss}
\setbox6\hbox{\hskip.6em\copy1\copy1$##1##1$}
\ifdim\MAXright<\wd6\MAXright\wd6\fi}
\def\norootupright##1{\hbox to0pt{\raise.45em\copy1\hskip-.36em\raise1.3\ht1%
\hbox{\copy4\hskip.3em$\ssize##1$}\hss}%
\setbox6\hbox{\hskip.6em\copy1\copy1$##1##1$}%
\ifdim\MAXright<\wd6\MAXright\wd6\fi}%
\def\rootdownright##1{\hbox to0pt{\raise-.5em\copy2\hskip-.25em\raise-1.35\ht1%
\hbox{\copy3\hskip.3em$\ssize##1$}\hss}\setbox6%
\hbox{\hskip.6em\copy2\copy2$##1##1$}%
\ifdim\MAXright<\wd6\MAXright\wd6\fi}%
\def\whiterootdownright##1{\hbox to0pt{\raise-.5em\copy2\hskip-.25em\raise-1.35\ht1
\hbox{\copy7\hskip.3em$\ssize##1$}\hss}\setbox6
\hbox{\hskip.6em\copy2\copy2$##1##1$}
\ifdim\MAXright<\wd6\MAXright\wd6\fi}
\def\rootdown##1{\hbox to0pt{\hskip-.05em\vrule height.25em depth.65em%
\hskip-.25em\raise-.95em\hbox{\copy3\hskip.3em$\ssize##1$}\hss}%
\setbox6\hbox{$##1$}%
\ifdim\MAXright<\wd6\MAXright\wd6\fi}%
\def\whiterootdown##1{\hbox to0pt{\hskip-.05em\vrule height.25em depth.65em
\hskip-.25em\raise-.95em\hbox{\copy7\hskip.3em$\ssize##1$}\hss}
\setbox6\hbox{$##1$}
\ifdim\MAXright<\wd6\MAXright\wd6\fi}
\def\dots{\hskip.5em\cdots\hskip.5em}}%
\def\enddynkin{\ifdim\MAXright>1em\hskip.5\MAXright\else\hskip.5em\fi\egroup}%
 
This macro is used for creation of Dynkin diagrams in AMSTeX.
The usage:
In math-mode, between \dynkin and \enddynkin,
\root#1 makes a node with label #1 over this node,
\rroot#1 and \droot do the same with labels on right and down,
\noroot and \norroot do the same with cross instead of bullet,
\link links two adjacent \root's, \llink#1 links two adjacent nodes
with double-line with #1 (probably '>' or '<') in the middle, \lllink is does
same with three lines, \dots does \cdots between the two adjacent nodes,
\rootdown#1 makes a linked node down labeled by #1, \rootupright#1 and
\rootdownright#1 are labeled linked nodes in the indicated directions 
which do not count for the jorizontal dimensions.

\sdynkin ... \endsdynkin works exactly as above but in smaller
version (eightpoint)

Examples:
$\dynkin \root{a_1}\link\root{a_2}\dots\root{a_{n-1}}\link\root{a_n}
\enddynkin$
or
$\dynkin \root{}\lllink>\root{}\enddynkin$
or
$\dynkin \root{a}\link\root{b}\rootupright{c}\rootdownright{d}\enddynkin$
or
$\dynkin \root{}\link\root{}\rootupright{}\link\root{}\enddynkin$.
or (white roots!)
{$\dynkin \whiteroot{0}\link\noroot{0}\dots\root{0}\link\whiteroot{0}\rootupright{0}\whiterootdownright{0}\enddynkin$}. 

\newcommand\C{{\mathbb C}}
\newcommand\E{{\mathbb E}}
\newcommand\HH{{\mathbb H}}
\newcommand\PP{{\mathbb P}}
 
 \newcommand\R{{\mathbb R}}

 \newcommand\bid{\mathbb{I}}
 \newcommand\Sp{{\mathbb S}}

 \newcommand\V{{\mathcal V}}

 \newcommand\g{{\mathfrak g}}
 \newcommand\h{{\mathfrak h}}
 \newcommand\p{{\mathfrak p}}
 \newcommand\q{{\mathfrak q}}
 \newcommand\gor{{\mathfrak r}}

 \newcommand\so{{\mathfrak{so}}}     
 \newcommand\slAlg{{\mathfrak{sl}}}

 \newcommand\al{\alpha}
 \newcommand\be{\beta}

 \newcommand\de{\delta}

 \newcommand\la{\lambda}

 \newcommand\La{\Lambda}

\newcommand\om{\omega}

\usepackage[all]{xy}

\newtheorem{thm}{Theorem}[section]

 \theoremstyle{definition}
 
 \theoremstyle{remark}
 
 \newtheorem*{ex}{Example}
 \numberwithin{equation}{section}
\newcommand\bV{{\mathbb V}}
\newcommand\bW{{\mathbb W}}
\newcommand\bU{{\mathbb U}}
\begin{document}

%
%
%
%
%
%
%
%
%
\title[A complex for the Dirac operator in several variables in dimension $4$]
{A complex for the Dirac operator in several variables in dimension $4$}
\author[L. Krump, V. Sou\v cek]{Luk\'a\v s Krump, Vladim\'{\i}r Sou\v cek}

\address{%
Mathematical Institute\\
Faculty of Mathematics and Physics\\
Charles University\\
Sokolovsk\'a 83\\
186 75 Praha\\
Czech Republic}

\email{krump@karlin.mff.cuni.cz, soucek@karlin.mff.cuni.cz}

\thanks{The work presented  here  was supported by the grant GA\v CR 20-11473S}

\begin{abstract}
The Penrose transform was used to construct a complex starting with the Dirac operator
in $k$ Clifford variables in dimension $2n$ in the stable range $n\geq k.$ 
In the paper, we consider  the same Penrose transform in the special case of dimension $4$  and for any number of variables (i.e., in the nonstable range).  In this case, we describe  explicitely the corresponding relative BGG complex and its direct image for cohomology with values in a suitable line bundle (in singular infinitesimal character).  We show that how to construct then a complex
starting with the Dirac operator in any number of variables. 

\end{abstract}

\maketitle
 
\section{Introduction.}

The motivation for the paper is coming from the theory of several Clifford variables.
The Dirac operator is broadly accepted as the appropriate generalization of the Cauchy-Riemann equations
to higher dimensions. Solutions of the Dirac equation have many similar properties to holomorphic functions.
The Clifford analysis describes properties of such solutions systematically (see \cite{BDS,DSS,Del}).

A study of solutions of the Dirac operator $D_k=(D^1,\ldots, D^k)$ in $k$ variables
is a higher-dimensional analogue of the theory of several complex variables. 
The problem studied in the paper is to construct  there  an analogue of the Dolbeault
complex for the case of the Dirac equation in several variables,
which starts with the operator $D.$ The answer is known in the case of 
in $k$ variables in dimension $2n$ in the stable range $n\geq k$
and it is visible that, unlike the Dolbeault case, there are second order operators
in the sequence.

The case of dimension 4 is quite exceptional. There are two strong reasons why it is so.
Firstly, the Euclidean space of dimension 4 can be identified with the field of quaternions. Quaternionic
analysis was first to be discovered and described appropriately (see \cite{Sud}). The case of several quaternionic variables was studied already in  \cite{Pert} and is now completely understood (see \cite{ALPS,Baston,BerSabStr,CSSS}). The second reason is deeper and is connected with symmetry of the corresponding system of equations.
In one variable, the symmetry group of the Dirac equation is the group of conformal symmetries of the sphere.
In more variables, the sphere is replaced by a homogeneous model of another parabolic geometry. It means that
the analogue of the sphere for more variables can be written as a homogeneous space $G/P,$ where $G$ is the group
of symmetries of the given equation and $P$ is its appropriate parabolic subgroup. In dimension 4, the homogeneous
model is the quaternionic projective space $\PP^k(\HH)$ and the group $G$ is $GL(k+1,\HH).$ In general case of dimension $2n$, the homogeneous model is a suitable isotropic Grassmannian $IG_{2n,k}(\R)$ and the group $G$ is $Spin(2n,k+2n;\R).$
In dimension 4, both symmetries are possible but the quaternionic one is broader and more effective.

A construction of  sequences of standard invariant differential operators on general manifolds with a given
parabolic structure was described in \cite{CSS} and simplified and extended in \cite{CD}. Invariant differential
operators can act, in general, only among sections of bundles induced by representations on the same
orbit of the affine action of the corresponding Weyl group.  
The construction mentioned above covers the case of regular orbits.
There is no such general construction available for orbits, which are singular.
The Dirac equation in several variables can be interpreted as an invariant differential operator
in a parabolic geometry (for details, see \cite{F})  but the corresponding orbit of the Weyl group is singular.

Resolutions starting with $D$ were already constructed (in some cases) by various methods.
One is an algebraic approach summarized in \cite{CSSS}, which covers quite a few cases. 
Another recent approach is based on a study of homomorphisms of (generalized) Verma modules, see \cite{F}.
Here we are starting an approach based on the Penrose transform on homogeneous spaces.
The method of construction of the complexes based on the Penrose transform was developed
by R. Baston and M. Eastwood (\cite{BE}). Inspiration for us came from the paper by R. Baston (\cite{Baston}) who used the Penrose transform as a tool in the case of quaternionic manifolds.

In the previous paper \cite{ks1}, we have already applied the Baston method based on the Penrose transform  in the case, where the conformal sphere was chosen to be the twistor space. It was shown there how to construct certain complexes
of invariant operators in singular character. In this case, the isotropy group was the smallest possible, with two elements.
The case of two Dirac operators is also in singular character but it belongs to an orbit with a bigger isotropic group having
four elements. 

A choice of a twistor space is given by a choice of a second parabolic subgroup in $G^c=Spin(2,k+2,\C),$ i.e. by  a choice
a second cross at the corresponding Dynkin diagram. The choice of the twistor space considered in \cite{ks1} covers only
orbits with the property that the isotropic subgroup has at most two elements. In this paper we show that another choice of the twistor
space makes it possible to cover more degenerate singular orbits. One of cases obtained here includes the complex starting
with the Dirac operator in two variables.

   Geometry of the Penrose transform on homogeneous spaces is based on a double fibration of suitable homogenous spaces and the transform translates objects on the twistor space (usually called the left bottom space) to the right bottom space.
   We shall use the transform for complex Lie groups, i.e. we shall consider 
   a complex simple Lie group $G^c$  and its parabolic subgroups. We shall formulate our problem in the real situation
and for proofs, we shall complexify our spaces and we shall use the machinery of the complex version of the Penrose transform.

The quaternionic resolution is (see ...) (here "$\implies$" denotes a second order operator):

$$
\bid\otimes\C_2\rightarrow\bid\otimes\C_{2k} \implies\bid\otimes\La^3\C_{2k} \rightarrow\C_2\otimes\La^4\C_{2k} \rightarrow$$
$$\rightarrow\odot^2\C_2\otimes\La^5\C_{2k} \rightarrow\dots\rightarrow\odot^{2k-3}\C_2\otimes\La^{2k}\C_{2k}\simeq\odot^{2k-3}\C_2$$

\section{The Penrose transform in principle.}

All necessary information concerning the Penrose transform can be found in the book \cite{BE}.
In the general setting, the geometry of the Penrose transform on homogeneous spaces
is given by a double fibration of homogeneous spaces
$$
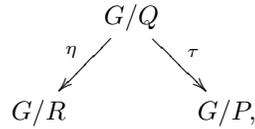
\xymatrix@C3mm@R7mm{
& G/Q \ar[dl]_{\eta} \ar[dr]^{\tau} \\
G/R && G/P,
} 
$$
where $P, Q, R$ are parabolic subgroups of a semisimple complex Lie group $G$ with $Q=P\cap R$.
 The main purpose of the paper is to construct certain complexes of $G$-invariant differential
 operators on the flag manifold $G/P$ for a suitable choice of the double fibration, so
 we shall restrict neglect some other features of the Penrose transform machinary.
 
 The construction has two steps. 
 Firstly,  for any irreducible $R$-module $\bV,$ there is a  locally exact complex of $G$-invariant
 differential operators on the flag manifold $G/Q$ called {\it the relative BGG resolution.}
This complex was introduced in \cite{BE} and the construction of relative BGG sequences of invariant operators was developed with all details and for (possibly curved) Cartan geometries
\cite{CS}.  An intuitive picture behind the relative BGG resolution is that classical BGG resolutions
along individual fibres of the map $\eta$  are glued together following the structure
of the fibration $\eta:G/Q\to G/R$ into the global complex of $G$-invariant differential operators  on $G/Q.$
An important advantage of the construction is that it applies both in regular and singular
infinitesimal character.

The second step of the contruction is to push down the complex given by the relative BGG resolution on $G/Q$ using the map $\tau $ to the flag manifold $G/P.$ Details of the procedure
in the case needed for our applications can be found in the paper \cite{Sal}.

For an explicit construction of the complexes
  the double fibration will be restricted
to open subsets.
We start with the big cell $U_p \subset G/P$   and we shall consider its $\tau$-preimage  $U_q=\tau^{-1}(U_p)\subset G/Q$ and the $\eta$-image $U_r=\eta(U_q)\subset G/R$ of $U_q$.
The double fibration has then the form
$$
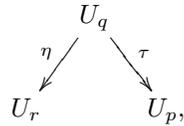
\xymatrix@C3mm@R7mm{
	& U_q\ar[dl]_{\eta} \ar[dr]^{\tau} \\
	U_r && U_p,
} 
$$
In this local situation, we can replace sections of vector bundles by maps with values in inducing modules.

\section{Notation.}

For the case we want to study, we consider the group $G=Spin(2k+4;\C)$.
The Lie algebra of $G$ is $\g=\so_{2k+4}\C$ and its Dynkin diagram with numbered roots is
$${\dynkin \root{1}\link\root{2}\dots\root{k-1}\link\root{k}\rootupright{k+1}\rootdownright{k+2}
\enddynkin}.$$

The parabolic subgroups $P,Q,R$ are given by the semisimple parts $\p^s, \q^s, \gor^s$ of Levi factors of the Lie algebras $\p, \q, \gor$. These can be given by crossing certain nodes, namely 

The standard  parabolic subalgebras $\p,\q,\gor$ are defined using their standard Dynkin diagram notation
$$\begin{aligned}
\p& = {\dynkin \root{}\link\root{}\dots\root{}\link\noroot{}\rootupright{}\rootdownright{}
\enddynkin} \\ 
\q&  = {\dynkin \root{}\link\root{}\dots\root{}\link\noroot{}\norootupright{} \rootdownright{} \enddynkin}\\ 
\gor&  = {\dynkin \root{}\link\root{}\dots\root{}\link\root{}\norootupright{}\rootdownright{}
\enddynkin}
\end{aligned}
$$
and they are sums of the Levi part and the nilpotent radical, the Levi part is the sum of its semisimple part and a commutative part:
$$
\begin{aligned}
\p=\p_0\oplus\p_+;\;\;\p_0=\p_0^{ss}\oplus \C;\;\;\p_0^{ss}\simeq \slAlg_k\C \times \so_4\C ;\\
\q=\q_0\oplus\q_+;\;\;\q_0=\q_0^{ss}\oplus \C_2;\;\q_0^{ss}\simeq\slAlg_k\C \times \slAlg_2\C ;\\
\gor=\gor_0\oplus\gor_+;\;\;\;\gor_0=\gor_0^{ss}\oplus \C;\;\;\;\gor_0^{ss}\simeq\slAlg_{k+2}\C.\;\;\;\;\;\;\;\;
\end{aligned}
$$
The parabolic subgroups $P,Q,R\subset G$ are fixed by the requirement that their Lie algebras
are equal to $\p,\q,\gor,$ respectively.
 
Note that while $P$ is determined by the symmetry of the given problem, different choices of $R$ lead
to different results and the choice of $R$ should be adapted to the  situation.   Finally, $Q=P\cap R$.

We suppose that  the Cartan subalgebra of $\h\subset \g$ is included in $ \p_0, \q_0,$  and $ \gor_0.$   Roots of $\g$ are 
$$
\pm\al_{ij}=\pm(e_i-e_j),\,\pm\be_{ij}=\pm(e_i+e_j), \text{ with } i,j=1,\ldots,k+2;\; i < j, 
$$
where $\{e_i\}$ is the standard basis of the dual $\h^*$.
Simple roots are $$\al_{12},\; \dots,\; \al_{k+1,k+2},\; \be_{k+1,k+2}.$$

For all of the three homogeneous spaces in the scheme, we shall consider representations of the corresponding parabolic subalgebras $\p, \q, \gor$.
Their irreducible representations are classified by irreducible representations of $\p_0, \q_0, \gor_0,$ which are in turn characterized by their highest weights $\la=[\la_1,\dots,\la_{k+2}]=\sum_{i=1}^{k+2}\la_i e_i$.
The action of root reflections on $\la$ (which will be used in the Penrose transform procedure) is given by permutations and sign changes of the coordinates of $\la$. 

Weights have the same number of coordinates for all three parabolic subalgebras.
To indicate which subalgebra we have in mind, we shall use the following convention due to M. Eastwood. 
In the Dynkin diagram notation, parabolic subalgebras are indicated by the appropriate crosses. 
In the notation $\la=[\la_1,\dots,\la_{k+2}]$, we shall use vertical lines after the $k$-th and/or $(k+1)$-st coordinate to indicate the crosses over the $k$-th and/or $(k+1)$-st node, respectively. While a weight $\la=[\la_1,\dots,\la_{k+2}]$ is $\g$-dominant if $\la_1\geq\dots\geq\la_k\geq\la_{k+1}\geq|\la_{k+2}|$, every cross in the Dynkin diagram cancels some of the inequalities. This is specified in the following table. 

\begin{center}
\begin{tabular}{|l|c|l|}
\hline
notation & dominant for & satisfies \\
\hline
$\la^p=[\la_1,\dots,\la_k|\la_{k+1},\la_{k+2}]$ & $\p^s$ & $\la_1\geq\dots\geq\la_k$ and $\la_{k+1}\geq|\la_{k+2}|$ \\
$\la^q=[\la_1,\dots,\la_k|\la_{k+1}|\la_{k+2}]$ & $\q^s$ & $\la_1\geq\dots\geq\la_k$ and $\la_{k+1}\geq-\la_{k+2}$ \\
$\la^r=[\la_1,\dots,\la_k,\la_{k+1}|\la_{k+2}]$ & $\gor^s$ & $\la_1\geq\dots\geq\la_k\geq\la_{k+1}\geq-\la_{k+2}$ \\
\hline
\end{tabular}
\end{center}

It is important to note that we follow the convention used in \cite{BE} and we
denote by $\E_{\la}$ the {\it dual} of the representation with the highest weight $\la,$
or equivalently, $\E_\la$ is the representation with the lowest weight $-\la.$
We shall also use the alternative notation for weights, where a weight $\la$ is given by the
coefficients in its decomposition into a linear combination of fundamental weights. Coefficients are then written over the corresponding nodes of the Dynkin diagrams.  
Our convention for fundamental weights is 
$$
\om_1=[1,0,\dots,0],\; 
\om_2=[1,1,0,\dots,0],\; \dots,\;
\om_k=[1,\dots,1,0,0], $$
$$ \om_{k+1}=\frac{1}{2}[1,\dots,1,-1],\;
\om_{k+2}=\frac{1}{2}[1,\dots,1,1].
$$
We also introduce the lowest positive weight $\de=\sum_{i=1}^{k+2}\om_i=[k+1,\dots,1,0]$ which will be used in the
affine action of the Weyl group.
 
\section{The relative BGG sequence.}

In this section, we follow the Penrose transform algorithm for our case (for details see \cite{BE}, chapter 4). 
The computation of the relative BGG sequence in a general situation of $k$ Dirac operators
in dimension $2n$ can be found in \cite{Sal2}. For $2n=4,$ the description simplifies substantially
and it is possible to describe the resolution very explicitely.
 
The first step is to determine the relative Hasse diagram which corresponds to the Hasse diagram of a typical fiber of the fibration $\eta,$ which is diffeomorphic with the Grassmannian of $2$-dimensional subspaces in $\C_{k+2}.$  In the Dynkin diagram
 notation it corresponds to
$$
	{\dynkin \root{}\link\root{}\dots\root{}\link\noroot{}\link\root{}\enddynkin} 
$$
(where the simple roots are $\al_{12},\;\dots,\;\al_{k,k+1},\;\be_{k+1,k+2}$). 
The BGG resolution is well known for this case.

For $k$ even, the Hasse diagram looks as follows (notation explained below): 


$$
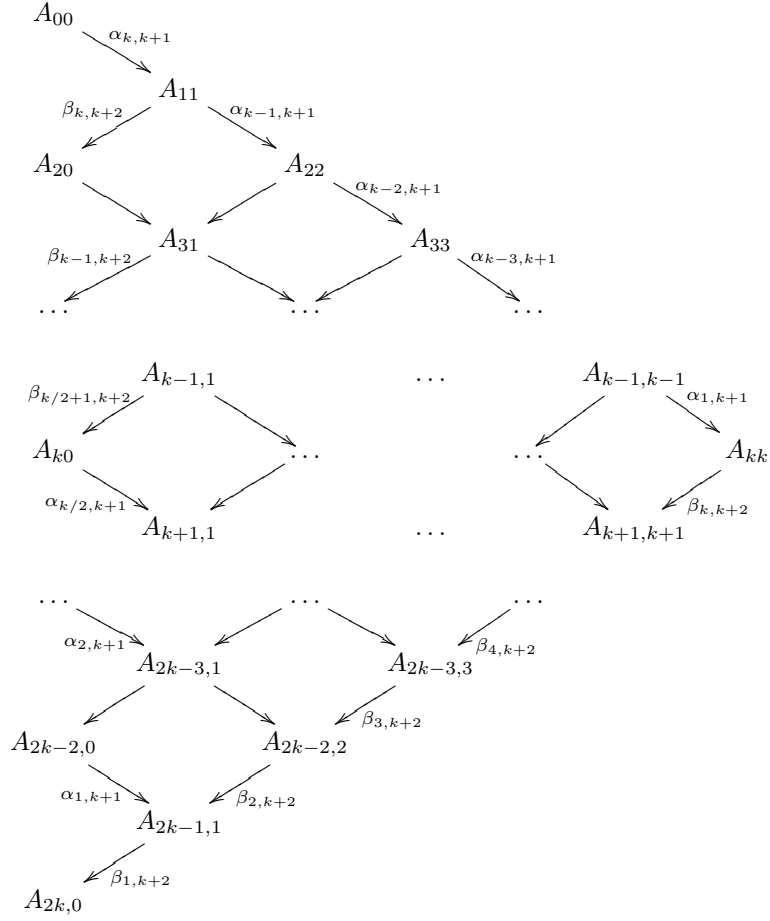
\xymatrix@C3mm@R5mm{
	A_{00} \ar[dr]^{\al_{k,k+1}} \\
	& A_{11} \ar[dl]_{\be_{k,k+2}} \ar[dr]^{\al_{k-1,k+1}}  \\
	A_{20} \ar[dr] && A_{22} \ar[dl] \ar[dr]^{\al_{k-2,k+1}}& \\
	& A_{31} \ar[dl]_{\be_{k-1,k+2}} \ar[dr] & & A_{33}\ar[dl] \ar[dr]^{\al_{k-3,k+1}}\\
	\dots & & \dots & & \dots\\
	& A_{k-1,1} \ar[dl]_{\be_{k/2+1,k+2}} \ar[dr] & & \dots  & & A_{k-1,k-1} \ar[dl] \ar[dr]^{\al_{1,k+1}}\\
	A_{k0} \ar[dr]_{\al_{k/2,k+1}} & & \dots \ar[dl]  & & \dots  \ar[dr] & & A_{kk}\ar[dl]^{\be_{k,k+2}}\\
	& A_{k+1,1}  & & \dots  & & A_{k+1,k+1} \\
	\dots \ar[dr]_{\al_{2,k+1}} & & \dots \ar[dl] \ar[dr] & & \dots \ar[dl]^{\be_{4,k+2}}\\
	& A_{2k-3,1} \ar[dl] \ar[dr] & & A_{2k-3,3}\ar[dl]^{\be_{3,k+2}}\\
	A_{2k-2,0} \ar[dr]_{\al_{1,k+1}} & & A_{2k-2,2} \ar[dl]^{\be_{2,k+2}}\\
	& A_{2k-1,1} \ar[dl]^{\be_{1,k+2}} \\
	A_{2k,0}
}
$$

For $k$ odd, there is a difference in the central three lines, which look like 

$$
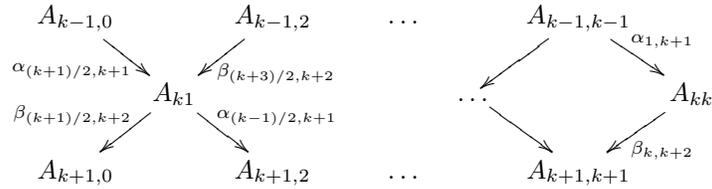
\xymatrix@C3mm@R5mm{
	A_{k-1,0}  \ar[dr]_{\al_{(k+1)/2,k+1}} & & A_{k-1,2} \ar[dl]^{\be_{(k+3)/2,k+2}}& & \dots  & & A_{k-1,k-1} \ar[dl] \ar[dr]^{\al_{1,k+1}}\\
	& A_{k1} \ar[dl]_{\be_{(k+1)/2,k+2}} \ar[dr]^{\al_{(k-1)/2,k+1}} & & & & \dots  \ar[dr] & & A_{kk}\ar[dl]^{\be_{k,k+2}}\\
	A_{k+1,0}  & & A_{k+1,2} & & \dots  & & A_{k+1,k+1}
}
$$

Note that roots assigned to the arrows are identical over parallel arrows, so that any (oriented) path between two vertices in this graph contains the same set of roots.

The second step is to turn this Hasse diagram into a BGG diagram in the following way: if $\la^q_0$ is a weight for the subalgebra $\q$, the $\q$-modules $A_{ij}$ are specified in the following way: let 
$$
\begin{aligned}
 	A_{00}&:=\E_{\la_0}, \\
	A_{11}&:=\E_{\mu}\text{ for }\mu=\la_0-\al_{k,k+1}, \\
	A_{20}&:=\E_{\nu}\text{ for }\nu=\la_0-\al_{k,k+1}-\be_{k,k+2} \\
	\text{ etc. } 
\end{aligned}
$$
More precisely, for every index $ij$ there is a unique set $R$ of roots, which are assigned to arrows between $A_{00}$ and $A_{ij}$ (as explained above), and then $A_{ij}=\E_\nu$ for $\nu=\la_0-\sum_{\al\in R}\al$.

The BGG diagram gives information on $\q$-modules, in which mappings in the corresponding relative BGG sequence take values. The weights in this diagram are all $\q$-dominant.

Let us define now
\begin{equation}\label{bV}
\bV_j=\sum_{i,0\leq 2i\leq j}A_{j,j-2i};\;j=0,1,\ldots,k;\;\;
\bV_{k+j}=\sum_{i,0\leq 2i\leq j}A_{k+j,k-j-2i};\;j=1,\ldots,k 
\end{equation}
and let $\V_j,j=0,\ldots 2k$  denote the sheaf of holomorphic sections of vector bundles
associated to modules $\bV_j.$ 

 For an application in Clifford analysis, we shall choose the initial weight as
$$ 
\la_0=
{\dynkin \root{0}\link\root{0}\dots\root{0}\link\root{0}\norootupright{-3}\rootdownright{0}\enddynkin}=\frac12[-3,\dots,-3|-3|3],
$$ 

\begin{thm}[Relative BGG resolution]\label{relBGG}
	Let us consider one-dimensional $Q$-module 
	$\E_{\lambda_0}$ for the weight
	 $\lambda_0=\frac12[-3,\dots,-3|-3|3].$
	Then there is an exact sequence
$$
\V_0\stackrel{D_1}{\longrightarrow}
\V_1\stackrel{D_2}{\longrightarrow}\V_2\ldots\V_{2k-1}\stackrel{D_{2k}}{\longrightarrow}\V_{2k}
$$
where $D_i$ are $G$-invariant differential operators of first order.	
\end{thm}

\vskip2mm\noindent{\bf Proof.}
The relative BGG resolution for the trivial representation $\E_{\lambda},\lambda=[0,\ldots,0]$
is just the de Rham resolution. The case considered in theorem is just the de Rham sequence
twisted by one-dimensional module $\E_{\lambda_0}.$ Details comments on the construction
of the resolution in the theorem can be found in \cite{Sal2}.

\hfill$\square$

\begin{ex}

As an illustration of the relative BGG sequence let us show the case $k=3$ (omitting the factor $\frac12$ in every weight):

$$
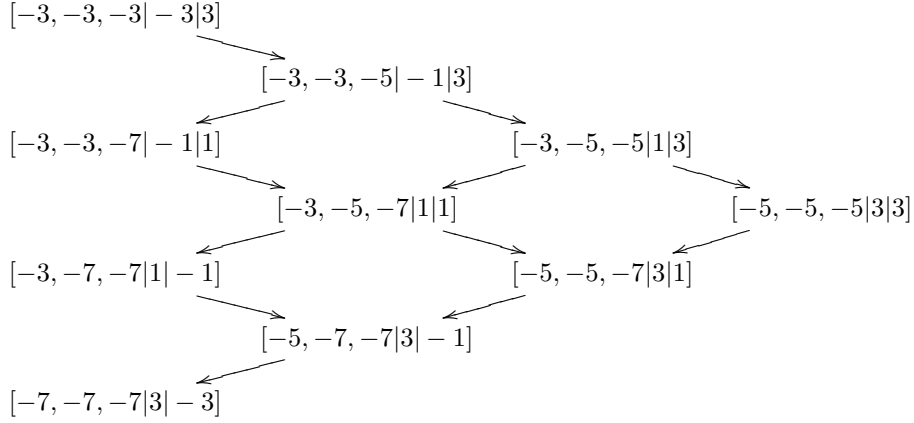
\xymatrix@C3mm@R3mm{
	[-3,-3,-3|-3|3] \ar[dr] \\
	& [-3,-3,-5|-1|3] \ar[dl] \ar[dr] \\
	[-3,-3,-7|-1|1] \ar[dr] && [-3,-5,-5|1|3]  \ar[dl] \ar[dr]& \\
	& [-3,-5,-7|1|1] \ar[dl] \ar[dr] & & [-5,-5,-5|3|3] \ar[dl] \\
	[-3,-7,-7|1|-1] \ar[dr] & & [-5,-5,-7|3|1] \ar[dl] \\
	& [-5,-7,-7|3|-1] \ar[dl] \\
	[-7,-7,-7|3|-3] 
}
$$
\end{ex}

Note that the chosen $\la_0$ is one of those for which the orbit belongs to a singular character. The first $k$ (in our example, the first $3$) coefficients indicate a $\slAlg_k\C$-dominant weight (since these are determined up to an additive constant), and the whole weight is obviously $\q$-dominant.

\section{Pushing down the relative BGG complex}

 \subsection{Calculation of direct images of weights}
Having a BGG sequence for $G/Q$, we can perform the third step: compute direct images of every individual weight. This is realised as the affine action of certain elements of the Weyl group (see again \cite{BE}, chapter ...). In our case, 
there is just one such element, namely the reflection in the direction of the root $\al_{k+1,k+2}$, let us call it $w$. The formula for the affine action is therefore $\la \mapsto \phi(\la)=w(\la+\de)-\de$. In practice, in the expression $$\la+\de=[\la_1+k+1,\dots,\la_k+2|\la_{k+1}+1|\la_{k+2}]$$ the last two coordinates are transposed (action of $w$) and $\de$ is subtracted. The resulting weight is regarded as a $\p$-weight.
Then we have the following possibilities: 
\begin{itemize}
\item either $\la$ is $\p$-dominant, then the direct image of $\la^q$ is $\la^p$ in degree $0$,
\item or $\phi(\la)$ is $\p$-dominant, then the direct image of $\la^q$ is $\phi(\la)^p$ in degree $1$, 
\item or $\la=\phi(\la)$, then there is no direct image.
\end{itemize}

Here is the example $k=3$, we omit the factor $\frac12$ in every weight again, $\emptyset$ stands for `no direct image', the subscript indicates the degree of the image:

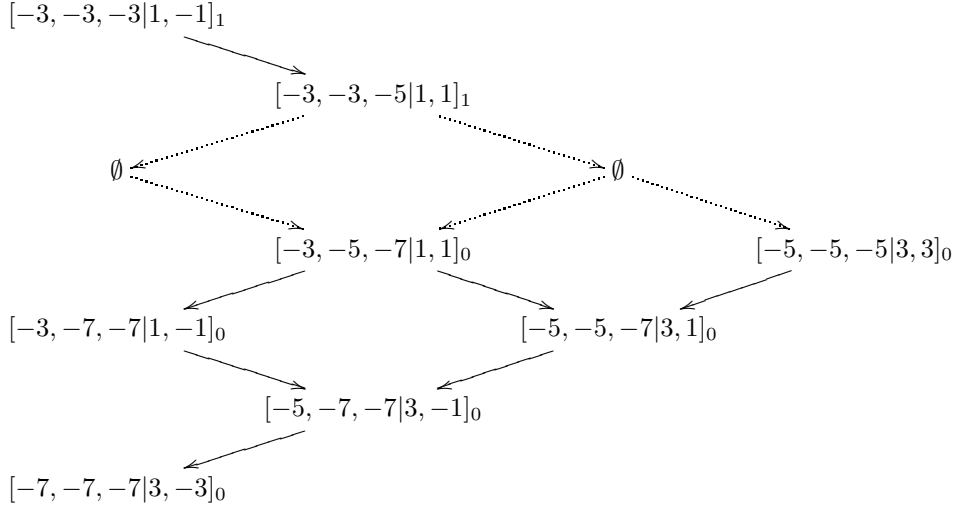
\begin{ex}

$$
\xymatrix@C3mm@R5mm{
[-3,-3,-3|1,-1]_1 \ar[dr] \\
& [-3,-3,-5|1,1]_1 \ar@.[dl] \ar@.[dr] \\
\emptyset \ar@.[dr] && \emptyset  \ar@.[dl] \ar@.[dr]& \\
& [-3,-5,-7|1,1]_0 \ar[dl] \ar[dr] & & [-5,-5,-5|3,3]_0 \ar[dl] \\
[-3,-7,-7|1,-1]_0 \ar[dr] & & [-5,-5,-7|3,1]_0 \ar[dl] \\
& [-5,-7,-7|3,-1]_0 \ar[dl] \\
[-7,-7,-7|3,-3]_0
}
$$
\end{ex}

Note that the $\slAlg_k\C$-part $\frac12[\la_1,\dots,\la_k|$ can also be expressed with nonnegative coefficients just adding a common constant. For instance, in the weight on the position $11$, we have $\frac12[-3,-3,-5|\simeq\frac12[2,2,0|=[1,1,0|$, which denotes $\C^3$, the basic representation of $\slAlg_3\C$. The $\so_4\C$-part in our example is $\frac12|1,1]$, i.e. the spinor representation $\Sp_+$. In such a way, it is easy to describe all the direct images for a general $k$.

Let us denote $\bU_{ij}=V_{ij}\otimes W_{ij}$ the $(\slAlg_k\C \times \so_4\C)$-module in the $i$-th row and $j$-th column of the sequence, where $i\in\{0,\dots,2k\}, j\in\{0,\dots,\min(i,2k-i)\}, i\neq2, i-j$ even. 

For the rows $i=0,1$ there are direct images in degree one. The module $U_{00}$ is obviously $\bid\otimes\Sp_- = \Sp_-$
($\bid$ is the trivial representation given by $[0,\dots,0|$), while $U_{11}$ is $\C^k\otimes\Sp_+$, therefore the initial arrow corresponds exactly to the Dirac operator $D_k$ (there is a unique possible operator -- see for instance \cite{F}).

In the the row $2$ there is no direct image, so we take $i\in\{3,\dots,2k\}$ in the sequel. 

For all $i\in\{3,\dots,2k\}, j\in\{0,\dots,\min(i,2k-i)\}, i-j$ even,
$V_{ij}$ is the (dual of) the $\slAlg_k\C$-module given by the highest weight $$[2,\dots,2,1,\dots,1,0,\dots,0]$$ with $(k-\frac{i+j}2)$ times $2$, $j$ times $1$ and $\frac{i-j}2$ times $0$ (together $k$ entries). $W_{ij}$ is the $\so_4\C$-module with the highest weight $(i-3)\om'_{k+1}+j\om'_{k+2}$, where $\om'_{k+1}=\frac12[1,-1], \om'_{k+2}=\frac12[1,1]$, the fundamental weights of $\so_4\C=\slAlg_2\C\times\slAlg_2\C$, are restrictions of $\om_{k+1}, \om_{k+2}$ to the last two coordinates. Note also that in the case of orthogonal algebras, modules and their duals coincide.

The simple or double arrows indicate whether the operators are of order one or two, respectively. The second order operators appear between lines $1$ and $3$, since there is a `order jump', ore more exactly this can be proved checking the conformal weights (see \cite{F}). All information obtained here can be summarised in a picture for a general even $k$:

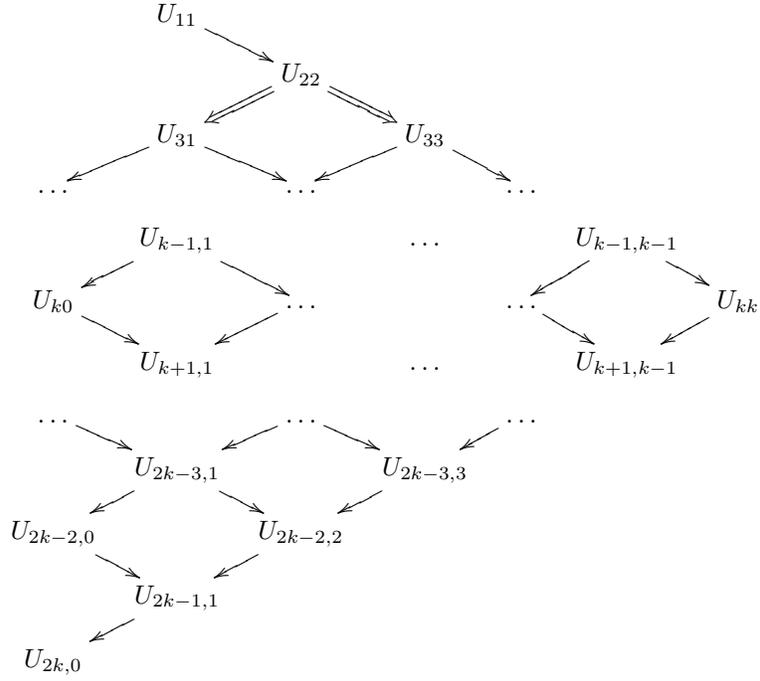
$$
\xymatrix@C3mm@R3mm{
& U_{11}\ar[dr]^{} \\
&& U_{22}\ar@2[dl]  \ar@2[dr]& \\
& U_{31} \ar[dl] \ar[dr] & & U_{33}\ar[dl] \ar[dr]\\
\dots & & \dots & & \dots\\
& U_{k-1,1} \ar[dl] \ar[dr] & & \dots  & & U_{k-1,k-1} \ar[dl] \ar[dr]\\
U_{k0} \ar[dr] & & \dots \ar[dl]  & & \dots  \ar[dr] & & U_{kk}\ar[dl]\\
& U_{k+1,1}  & & \dots  & & U_{k+1,k-1} \\
\dots \ar[dr] & & \dots \ar[dl] \ar[dr] & & \dots \ar[dl]\\
& U_{2k-3,1} \ar[dl] \ar[dr] & & U_{2k-3,3}\ar[dl]\\
U_{2k-2,0} \ar[dr] & & U_{2k-2,2} \ar[dl]\\
& U_{2k-1,1} \ar[dl] \\
U_{2k,0}
}
$$

And for a general odd $k$:

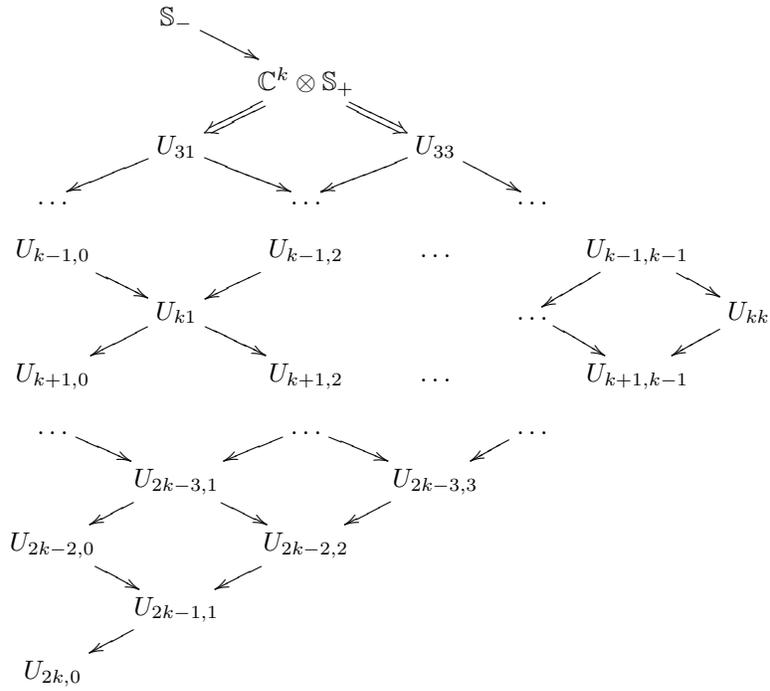
$$
\xymatrix@C3mm@R3mm{
& \Sp_- \ar[dr]^{} \\
&& \C^k\otimes\Sp_+ \ar@2[dl]  \ar@2[dr]& \\
& U_{31} \ar[dl] \ar[dr] & & U_{33}\ar[dl] \ar[dr]\\
\dots & & \dots & & \dots\\
U_{k-1,0}  \ar[dr] & & U_{k-1,2} \ar[dl]& \dots  & & U_{k-1,k-1} \ar[dl] \ar[dr]\\
& U_{k1} \ar[dl] \ar[dr] & & & \dots  \ar[dr] & & U_{kk}\ar[dl] \\
U_{k+1,0}  & & U_{k+1,2} & \dots  & & U_{k+1,k-1}\\
\dots \ar[dr] & & \dots \ar[dl] \ar[dr] & & \dots \ar[dl]\\
& U_{2k-3,1} \ar[dl] \ar[dr] & & U_{2k-3,3}\ar[dl]\\
U_{2k-2,0} \ar[dr] & & U_{2k-2,2} \ar[dl]\\
& U_{2k-1,1} \ar[dl] \\
U_{2k,0}
}
$$

\begin{ex}
For instance, for $k=5$:

$$
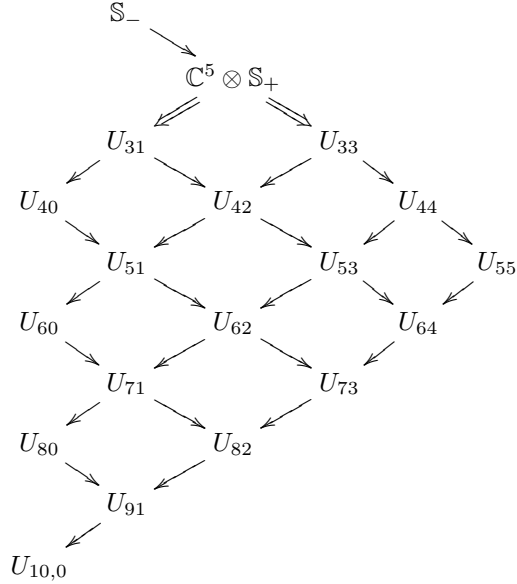
\xymatrix@C3mm@R3mm{
& \Sp_- \ar[dr]^{} \\
&& \C^5\otimes\Sp_+ \ar@2[dl]  \ar@2[dr]& \\
& U_{31} \ar[dl] \ar[dr] & & U_{33}\ar[dl] \ar[dr]\\
U_{40}\ar[dr] & & U_{42} \ar[dl] \ar[dr] & & U_{44}\ar[dl] \ar[dr]\\
& U_{51} \ar[dl] \ar[dr] & & U_{53} \ar[dl] \ar[dr] & & U_{55}\ar[dl]\\
U_{60} \ar[dr] & & U_{62} \ar[dl] \ar[dr] & & U_{64} \ar[dl]\\
& U_{71} \ar[dl] \ar[dr] & & U_{73} \ar[dl]\\
U_{80} \ar[dr] & & U_{82} \ar[dl]\\
& U_{91} \ar[dl]\\
U_{10,0}
}
$$

\end{ex}

\subsection{Restriction to the big cell $\g_{-}$}

To come closer to the standard setting of Clifford analysis, we shall use the 'big cell' in the flag manifold $G/P$
to make all vector bundles trivial and to identify the space of sections over the big cell
with spaces of holomorphic maps on the big cell with values in inducing representations.

The parabolic subalgebra $\p$ is related to the $|2|$-grading 
$\g=\g_{-2}\oplus\g_{-1}\oplus\g_{0}\oplus\g_{1}\oplus\g_{2}
$
by $\p_+=\g_1\oplus\g_2.$ Let $\g_-=\g_{-2}\oplus\g_{-1}$ and $G_-=\exp(\g_1).$

The projection $\pi:G\to M=G/P$ is injective on $G_-$ and its inverse $q$ is a section for $\pi$
restricted to $\pi(G_-).$ If $V$ is a homogeneous bundle  on $M$ induced by $P$-module $\bV,$
then the space of holomorphic sections $\Gamma(G_-,V)$ will be identified with the space
of holomorphic maps from $\g_-$ to $\bV.$ For brevity, we shall denote it just by $\bV.$

Then we get the following direct image of the relative BGG resolution down on $\g_{-}.$

\begin{thm}\label{kDirac}
	Let us consider the following set of $P$-modules:
	
	$$\bU_1=\frac{1}{2}[-3,\ldots,-3|1,-1],\; \bU_2=\frac{1}{2}[-3,\ldots,-3,-5|1,1];$$	
	Highest weights for $\q$-modules $\bV_{i},i=3,\ldots,2k$ defined in (\ref{bV})
	are $\p$-dominant and we consider them as $P$-modules.
	
	Then the direct image of the relative BGG resolution from Th.\ref{relBGG}
	is a complex
		\begin{equation}
	\mathcal{U}_1\stackrel{D'_1}{\longrightarrow}
	\mathcal{U}_2\stackrel{D'_2}{\longrightarrow}\V_2\stackrel{D'_3}{\longrightarrow}\V_3
	\ldots\V_{2k-1}\stackrel{D'_{2k}}{\longrightarrow}\V_{2k}
	\end{equation}
	The operator $D'_2$ is second order and all other operators are first order.
\end{thm}

\vskip2mm\noindent {\bf Proof.}
Starting from operator $D'_3$ on, it follows immediately from the fact that the relative
BGG resolution is a complex, because the all direct images of values are in degree zero.
To check that the composition of two following operators vanishes for the first three operators
follows easily from the definition of higher differential in the corresponding double complex.

\hfill $ \square$

  \section{Reality condition}
  Up to now, we have used the holomorphic version of the Penrose transform
  as formulated in \cite{BE}. For applications to Clifford analysis, we want to
  restrict now to a real version.
  
  So let now start with a real group $G=Spin(k,k+4)$ and with its parabolic subgroup $P$
  chosen in such a way that the simple part of the Levi subgroup is equal to $GL(k,\mathbb{R})\times
  Spin(4).$

Then $G$ is a real form of the group $G^\C=Spin(2k+4,\C) $ and we can choose the parabolic
subgroup $P^\C\subset G^\C$ such that
the flag manifold $G^\C/P^\C$ is a complexification of the Grassmannian $G/P$ of maximal isotropic subspaces of dimension $k$ in $\mathbb{R}^{k,k+4}.$

Suppose now that $\bV$ and $\bW$ are complex irreducible modules for $P^\C.$

As described in detail in \cite{Sal}, Sect. 2.6, $G^\C$-invariant operators from
the associated bundle $V$ to the associated bundle $W$ on $g^\C/P^\C$ correspond to
$G$-invariant operators from the associated bundle $V$ to the associated bundle $W$
on $G/P.$

So we get a complex of $G$ invariant operators induced from the complex constructed
in Theorem \ref{kDirac} on $G^\C/P^\C$ to the homogeneous space $G/P.$
 
\section{Restriction of the complex to $\g_{-1}$}

 The final step we need to do is to reduce the complex in Theorem \ref{kDirac}
 to spaces of maps defined on the big cell $\g_-$ to
 maps defined only on $G_{-1}.$  This is treated in detail in \cite{Sal}, Prop. 2.2.
 
 Let $G=Spin_{k,k+4}$ and let
 $\g_i, i=-2,-1,0,1,2$ gives the (real) $|2|$-grading described in \cite{Sal}.
 Then the space $M(4,k,\mathbb{R})\simeq\g_{-1}$ is maped by exponential onto 
 the big cell $G_-.$ Denote by $q$ the projection of $G_-$ to the space
 $U=G_{-2}/G_-$ of right cosets.
 
 Then we have
 \begin{thm}[Prop. 2.2 in\cite{Sal}]
 	
 	Let $\bV$ and $\bW$ are vector spaces and let
 	$$
 	D:\mathcal{C}^\infty(G_-,\bV)\to\mathcal{C}^\infty(G_-,\bW)
 	$$
 	be a linear $G_-$ invariant differential operator from $\mathcal{U}_r(\g_-).$
 	
 	Then there is a unique linear operator 
 	$$\underline{D}:\mathcal{C}^\infty(U,\bV)\to\mathcal{C}^\infty(U,\bW)$$
 	such that the diagram
 	$$
 	\begin{CD}
 		 \mathcal{C}^\infty(G_-,\bV)@>D>> \mathcal{C}^\infty(G_-,\bW)\\
 		 @AAq^*A @AAq^*A\\
 		\mathcal{C}^\infty(U,\bV) @>\underline{D}>> \mathcal{C}^\infty(U,\bW)
 	\end{CD}
 	$$
 	commutes. Moreover, $\underline{D}$ is a homogeneous constant coefficient operator of order $r.$
 	
 \end{thm}

As a corollary, we get immediately that if upstairs we have a trivial composition $D_1\circ D_2,$
then also $\underline{D}_1\circ\underline{D}_2$ is trivial map.

Hence the complex in Theorem \ref{kDirac} induces a complex of maps defined
on the space $M(4,k,\mathbb{R}).$ It is immediate to check that the first map $\underline{D}_1$
is the map given by the Dirac operator in $k$ variables in dimension $4.$

\end{document}